\documentstyle[12pt,amssymb]{article}

\newcommand{\Label}{\label} 
 
\begin{document} 
 
\baselineskip18pt 
 
\newtheorem{definition}{Definition $\!\!$}[section] 
\newtheorem{prop}[definition]{Proposition $\!\!$} 
\newtheorem{lem}[definition]{Lemma $\!\!$} 
\newtheorem{corollary}[definition]{Corollary $\!\!$} 
\newtheorem{theorem}[definition]{Theorem $\!\!$} 
\newtheorem{example}[definition]{\it Example $\!\!$} 
\newtheorem{remark}[definition]{Remark $\!\!$} 
 
\newcommand{\nc}[2]{\newcommand{#1}{#2}} 
\newcommand{\rnc}[2]{\renewcommand{#1}{#2}} 
 
\nc{\Section}{\setcounter{definition}{0}\section} 
\renewcommand{\theequation}{\thesection.\arabic{equation}} 
\newcounter{c} 
\renewcommand{\[}{\setcounter{c}{1}$$} 
\newcommand{\etyk}[1]{\vspace{-7.4mm}$$\begin{equation}\Label{#1} 
\addtocounter{c}{1}} 
\renewcommand{\]}{\ifnum \value{c}=1 $$\else \end{equation}\fi} 
 
%%%%%%% Piotr's macro %%%%%%%%%%%%%%%%%%%%%%%%%%%%%%%%%%%%%%%%%%%%%%%% 
\nc{\bpr}{\begin{prop}} 
\nc{\bth}{\begin{theorem}} 
\nc{\ble}{\begin{lem}} 
\nc{\bco}{\begin{corollary}} 
\nc{\bre}{\begin{remark}} 
\nc{\bex}{\begin{example}} 
\nc{\bde}{\begin{definition}} 
\nc{\ede}{\end{definition}} 
\nc{\epr}{\end{prop}} 
\nc{\ethe}{\end{theorem}} 
\nc{\ele}{\end{lem}} 
\nc{\eco}{\end{corollary}} 
\nc{\ere}{\hfill\mbox{$\Diamond$}\end{remark}} 
\nc{\eex}{\end{example}} 
\nc{\epf}{\hfill\mbox{$\Box$}} 
\nc{\ot}{\otimes} 
\nc{\bsb}{\begin{Sb}} 
\nc{\esb}{\end{Sb}} 
\nc{\ct}{\mbox{${\cal T}$}} 
\nc{\ctb}{\mbox{${\cal T}\sb B$}} 
\nc{\bcd}{\[\begin{CD}} 
\nc{\ecd}{\end{CD}\]} 
\nc{\ba}{\begin{array}} 
\nc{\ea}{\end{array}} 
\nc{\bea}{\begin{eqnarray}} 
\nc{\eea}{\end{eqnarray}} 
\nc{\be}{\begin{enumerate}} 
\nc{\ee}{\end{enumerate}} 
\nc{\beq}{\begin{equation}} 
\nc{\eeq}{\end{equation}} 
\nc{\bi}{\begin{itemize}} 
\nc{\ei}{\end{itemize}} 
\nc{\kr}{\mbox{Ker}} 
\nc{\te}{\!\ot\!} 
\nc{\pf}{\mbox{$P\!\sb F$}} 
\nc{\pn}{\mbox{$P\!\sb\nu$}} 
\nc{\bmlp}{\mbox{\boldmath$\left(\right.$}} 
\nc{\bmrp}{\mbox{\boldmath$\left.\right)$}} 
\rnc{\phi}{\mbox{$\varphi$}} 
\nc{\LAblp}{\mbox{\LARGE\boldmath$($}} 
\nc{\LAbrp}{\mbox{\LARGE\boldmath$)$}} 
\nc{\Lblp}{\mbox{\Large\boldmath$($}} 
\nc{\Lbrp}{\mbox{\Large\boldmath$)$}} 
\nc{\lblp}{\mbox{\large\boldmath$($}} 
\nc{\lbrp}{\mbox{\large\boldmath$)$}} 
\nc{\blp}{\mbox{\boldmath$($}} 
\nc{\brp}{\mbox{\boldmath$)$}} 
\nc{\LAlp}{\mbox{\LARGE $($}} 
\nc{\LArp}{\mbox{\LARGE $)$}} 
\nc{\Llp}{\mbox{\Large $($}} 
\nc{\Lrp}{\mbox{\Large $)$}} 
\nc{\llp}{\mbox{\large $($}} 
\nc{\lrp}{\mbox{\large $)$}} 
\nc{\lbc}{\mbox{\Large\boldmath$,$}} 
\nc{\lc}{\mbox{\Large$,$}} 
\nc{\Lall}{\mbox{\Large$\forall$}} 
\nc{\bc}{\mbox{\boldmath$,$}} 
\rnc{\epsilon}{\varepsilon} 
\rnc{\ker}{\mbox{\em Ker}} 
\nc{\ra}{\rightarrow} 
\nc{\ci}{\circ} 
\nc{\cc}{\!\ci\!} 
\nc{\T}{\mbox{\sf T}} 
\nc{\can}{\mbox{\em\sf T}\!\sb R} 
\nc{\cnl}{$\mbox{\sf T}\!\sb R$} 
\nc{\lra}{\longrightarrow} 
\nc{\M}{\mbox{Map}} 
%\rnc{\to}{\mapsto} 
\nc{\imp}{\Rightarrow} 
\rnc{\iff}{\Leftrightarrow} 
\nc{\bmq}{\cite{bmq}} 
\nc{\ob}{\mbox{$\Omega\sp{1}\! (\! B)$}} 
\nc{\op}{\mbox{$\Omega\sp{1}\! (\! P)$}} 
\nc{\oa}{\mbox{$\Omega\sp{1}\! (\! A)$}} 
\nc{\inc}{\mbox{$\,\subseteq\;$}} 
\nc{\de}{\mbox{$\Delta$}} 
\nc{\spp}{\mbox{${\cal S}{\cal P}(P)$}} 
\nc{\dr}{\mbox{$\Delta_{R}$}} 
\nc{\dsr}{\mbox{$\Delta_{\cal R}$}} 
\nc{\m}{\mbox{m}} 
\nc{\0}{\sb{(0)}} 
\nc{\1}{\sb{(1)}} 
\nc{\2}{\sb{(2)}} 
\nc{\3}{\sb{(3)}} 
\nc{\4}{\sb{(4)}} 
\nc{\5}{\sb{(5)}} 
\nc{\6}{\sb{(6)}} 
\nc{\7}{\sb{(7)}} 
\nc{\hsp}{\hspace*} 
\nc{\nin}{\mbox{$n\in\{ 0\}\!\cup\!{\Bbb N}$}} 
\nc{\al}{\mbox{$\alpha$}} 
\nc{\bet}{\mbox{$\beta$}} 
\nc{\ha}{\mbox{$\alpha$}} 
\nc{\hb}{\mbox{$\beta$}} 
\nc{\hg}{\mbox{$\gamma$}} 
\nc{\hd}{\mbox{$\delta$}} 
\nc{\he}{\mbox{$\varepsilon$}} 
\nc{\hz}{\mbox{$\zeta$}} 
\nc{\hs}{\mbox{$\sigma$}} 
\nc{\hk}{\mbox{$\kappa$}} 
\nc{\hm}{\mbox{$\mu$}} 
\nc{\hn}{\mbox{$\nu$}} 
\nc{\la}{\mbox{$\lambda$}} 
\nc{\hl}{\mbox{$\lambda$}} 
\nc{\hG}{\mbox{$\Gamma$}} 
\nc{\hD}{\mbox{$\Delta$}} 
\nc{\th}{\mbox{$\theta$}} 
\nc{\Th}{\mbox{$\Theta$}} 
\nc{\ho}{\mbox{$\omega$}} 
\nc{\hO}{\mbox{$\Omega$}} 
\nc{\hp}{\mbox{$\pi$}} 
\nc{\hP}{\mbox{$\Pi$}} 
\nc{\bpf}{{\sl Proof.~~}} 
\nc{\slq}{\mbox{$A(SL\sb q(2))$}} 
\nc{\fr}{\mbox{$Fr\llp A(SL(2,\IC))\lrp$}} 
\nc{\slc}{\mbox{$A(SL(2,\IC))$}} 
\nc{\af}{\mbox{$A(F)$}} 
\rnc{\widetilde}{\tilde} 
\nc{\qdt}{quantum double torus} 
\nc{\aqdt}{\mbox{$A(DT^2_q)$}} 
\nc{\dtq}{\mbox{$DT^2_q$}} 
\nc{\uc}{\mbox{$U(2)$}} 
\nc{\uq}{\mbox{$U_{q^{-1},q}(2)$}} 
\rnc{\subset}{\inc}

\def\esl{{\mbox{$E\sb{\frak s\frak l (2,{\Bbb C})}$}}} 
\def\esu{{\mbox{$E\sb{\frak s\frak u(2)}$}}} 
\def\zf{{\mbox{${\Bbb Z}\sb 4$}}} 
\def\zt{{\mbox{$2{\Bbb Z}\sb 2$}}} 
\def\ox{{\mbox{$\Omega\sp 1\sb{\frak M}X$}}} 
\def\oxh{{\mbox{$\Omega\sp 1\sb{\frak M-hor}X$}}} 
\def\oxs{{\mbox{$\Omega\sp 1\sb{\frak M-shor}X$}}} 
\def\Fr{\mbox{Fr}} 
\def\gal{-Galois extension} 
\def\hge{Hopf-Galois extension} 
\def\cge{coalgebra-Galois extension} 
\def\pge{$\psi$-Galois extension} 
\def\ta{\tilde a} 
\def\tb{\tilde b} 
\def\tc{\tilde c} 
\def\td{\tilde d} 
\def\st{\stackrel} 
 
%%%%%%%%%%%%%%%%% macro 2 %%%%%%%%%%%%%%%%%%%%%%%%%%%%%%%%%%%%%%%%%%%%%%%%% 
\newcommand{\Sp}{{\rm Sp}\,} 
\newcommand{\Mor}{\mbox{$\rm Mor$}} 
\newcommand{\skrA}{{\cal A}} 
\newcommand{\Phase}{\mbox{$\rm Phase\,$}} 
\newcommand{\id}{{\rm id}} 
\newcommand{\diag}{{\rm diag}} 
\newcommand{\inv}{{\rm inv}} 
\newcommand{\ad}{{\rm ad}} 
\newcommand{\poi}{{\rm pt}} 
\newcommand{\Dim}{{\rm dim}\,} 
\newcommand{\Ker}{{\rm ker}\,} 
\newcommand{\Mat}{{\rm Mat}\,} 
\newcommand{\Rep}{{\rm Rep}\,} 
\newcommand{\Fun}{{\rm Fun}\,} 
\newcommand{\Tr}{{\rm Tr}\,} 
\newcommand{\supp}{\mbox{$\rm supp$}} 
\newcommand{\half}{\frac{1}{2}} 
\newcommand{\skrF}{{A}} 
\newcommand{\skrD}{{\cal D}} 
\newcommand{\skrC}{{\cal C}} 
\newcommand{\ttimes}{\mbox{$\hspace{.5mm}\bigcirc\hspace{-4.9mm} 
\perp\hspace{1mm}$}} 
\newcommand{\Ttimes}{\mbox{$\hspace{.5mm}\bigcirc\hspace{-3.7mm} 
\raisebox{-.7mm}{$\top$}\hspace{1mm}$}} 
\newcommand{\bbr}{{\bf R}} 
\newcommand{\bbz}{{\bf Z}} 
\newcommand{\Ci}{C_{\infty}} 
\newcommand{\Cb}{C_{b}} 
\newcommand{\fa}{\forall} 
\newcommand{\rrr}{right regular representation} 
\newcommand{\wrt}{with respect to} 
\newcommand{\qg}{quantum group} 
\newcommand{\qgs}{quantum groups} 
\newcommand{\cs}{classical space} 
\newcommand{\qs}{quantum space} 
\newcommand{\po}{Pontryagin} 
\newcommand{\ch}{character} 
\newcommand{\chs}{characters} 
 
\def\inbar{\,\vrule height1.5ex width.4pt depth0pt} 
\def\IC{{\Bbb C}} 
\def\IZ{{\Bbb Z}} 
\def\IN{{\Bbb N}} 
\def\otc{\otimes_{\IC}} 
\def\ra{\rightarrow} 
\def\ota{\otimes_ A} 
\def\otza{\otimes_{ Z(A)}} 
\def\otc{\otimes_{\IC}} 
\def\h{\rho} 
\def\x{\zeta} 
\def\th{\theta} 
\def\s{\sigma} 
\def\t{\tau} 
 
%%%%%%%%%%%%%%%%%%% Tom's MACRO %%%%%%%%%%%%%%%%%%%%%%%%%%%%%%%% 
\def\sw#1{{\sb{(#1)}}} 
\def\sco#1{{\sp{(\bar #1)}}} % for coactions -- you can redefine 
% note its hard to read sometimes when can be confused with coproduct. 
\def\su#1{{\sp{(#1)}}} % for action of transl map and other situations 
\def\d{{\rm d}} 
\def\proof{{\sl Proof.~~}} 
\def\endproof{\hbox{$\sqcup$}\llap{\hbox{$\sqcap$}}\medskip} 
\def\tens{\mathop{\otimes}} 
\def\CC{{\cal C}} 
\def\CL{{\Lambda}} 
\def\CM{{\cal M}} % {\rm od}} 
\def\CN{{\cal N}} %  ideal of calculus on P -- you can redefine 
\def\o{{}_{(1)}} 
\def\t{{}_{(2)}} 
\def\Bo{{}_{\und{(1)}}} 
\def\Bt{{}_{\und{(2)}}} 
\def\th{{}_{(3)}} 
\def\Bth{{}_{\und{(3)}}} 
\def\<{{\langle}} 
\def\>{{\rangle}} 
\def\und#1{{\underline{#1}}} 
\def\la{{\triangleright}} % -- left, right actions 
\def\id{{\rm id}} % roman id map 
\def\eps{\epsilon} 
\def\span{{\rm span}} 
\def\q2{{q^{-2}}} 
\def\bicross{{\blacktriangleright\!\!\!\triangleleft}} 
\def\note#1{{}} 
\def\eqn#1#2{\begin{equation}#2\label{#1}\end{equation}} 
\def\qbinom#1#2#3{\left(\begin{array}{c}#1\\#2\end{array}\right)\sb#3} 
\def\Z{{\Bbb Z}} 
\def\can{{\rm can}} 
\def\note#1{} 
%%%%%%%%%%%%%%%%%%%%%%%%%%%%%%%%%%%%%%%%%%%%%%%%%%%%%%%%%%%%%%%%%%%%%%%%%%%%% 
~{\ }\qquad\hskip 3.8in \hskip -0.62802pt \phantom{DAMTP/97-76}\vspace{1in} 
\begin{center} 
{\large\bf COALGEBRA EXTENSIONS AND ALGEBRA 
COEXTENSIONS OF GALOIS TYPE} 
\vspace{12pt}\\ 
%%%%%%%%%%%%%%%%%%%%%%%%%%%%%%%%%%%%%%%%%%%%%%%%%%%%%
{\sc Tomasz Brzezi\'nski}\footnote{ 
Lloyd's fellow. On~leave from:  
Department of Theoretical Physics, University of \L\'od\'z, 
Pomorska 149/153, 90--236 \L\'od\'z, Poland. 
{\sc e-mail}: {\sc tb10@york.ac.uk}}\\
{Department of Mathematics,  University of York,\\
 Heslington, York YO1 5DD, U.K.}\vspace{5mm}\\ 
%%%%%%%%%%%%%%%%%%%%%%%%%%%%%%%%%%%%%%%%%%%%%%%%%%%
{\sc Piotr M.~Hajac}\footnote{ 
On~leave from: 
Department of Mathematical Methods in Physics, 
Warsaw University, ul.~Ho\.{z}a 74, Warsaw, \mbox{00--682~Poland}. 
 http://info.fuw.edu.pl/KMMF/ludzie$\underline{~~}$ang.html\\
{\sc e-mail}: {\sc hajac@sissa.it}} \\ 
{International School for Advanced Studies,\\ 
Via Beirut 2--4, 34013 Trieste, Italy.} 
\end{center} 
%\vspace{24pt} 

\begin{abstract} 
The notion of a coalgebra-Galois extension is defined as a natural 
generalisation of a Hopf-Galois extension. It is shown that any 
coalgebra-Galois extension induces a unique entwining map 
$\psi$ compatible with the right coaction. 
For the dual notion of an algebra-Galois coextension 
it is also proven that there always exists a unique entwining structure 
compatible with the right action. 
\end{abstract} 
 
%{\ }\qquad\qquad \hskip 4.3in DAMTP/97- 
%~{\ }\qquad\qquad \hskip 4.3in q-alg/97 
%~\vspace{.2in} 

\section{Introduction}\parindent5mm 
 
Hopf-Galois extensions  
can be viewed as  non-commutative torsors or principal 
bundles with universal differential structure. From the 
latter point of view a quantum group gauge theory 
was introduced in \cite{BrzMa:gau} and developed in~\cite{h,BrzMa:dif}.   
It turns out that to develop gauge 
theory on quantum homogeneous spaces (e.g., the family of Podle\'s 
quantum spheres) or with braided groups (see~\cite{Maj:adv} 
 for a recent review), one needs to consider coalgebra bundles.  
Such a generalisation of gauge theory was 
proposed recently in \cite{BrzMa:coa}. It introduces the notions 
of an entwining structure with an entwining map $\psi$, and  a  
$\psi$-principal coalgebra bundle. The latter can be viewed 
as a generalisation of a Hopf-Galois extension. 
 
In the present paper 
we introduce the notion of a coalgebra-Galois extension 
for  arbitrary coalgebras as a  natural 
generalisation of a Hopf-Galois extension. It is obtained by giving up
the condition that the coaction be 
%left-linear over the coinvariants rather than 
an algebra map. Our main result is that such defined coalgebra-Galois 
extension always induces a unique compatible entwining map $\psi$.  
We dualise coalgebra-Galois extensions and thus introduce the 
notion of an algebra-Galois coextension. We then prove the dual version
of the main result.
\note{
We show that for such a 
coextension there exists an entwining map $\psi$ making 
it a $\psi$-Galois coextension (called dual $\psi$-principal coalgebra bundle  
in~\cite{BrzMa:coa}). 
}
Finally, we prove that if two quotient coalgebras 
$C/I_1$ and $C/I_2$ cogenerate (in the sense of Definition~\ref{codef})  
the coalgebra $C$, then  
the coinvariants of $C$ are the same as the intersection of the coinvariants 
of $C/I_1$ with the coinvariants of~$C/I_2\,$. 
\note{in the 
context of gauge theory of quantum homogeneous spaces is, in fact, a 
very natural generalisation of Hopf-Galois extensions. More 
specifically, we show that the theory of $\psi$-principal bundles 
extends the Hopf-Galois extension theory (see \cite{Mon:hop} for a 
review) to the case in which the coaction is a left linear over the 
coinvariants rather than an algebra map.} 
 
{\em Notation.} Here and below 
$k$ denotes a field. All algebras are over $k$, 
 associative and unital with the unit denoted by~1. All algebra homomorphisms 
are assumed to be unital. We use the standard 
algebra and coalgebra notation, i.e., $\Delta$ is a coproduct, $m$ is a 
product, $\eps$ is a counit, etc. The identity map from the space $V$ to 
itself 
is also denoted by $V$. The unadorned tensor product stands for the tensor 
product over $k$. To abbreviate notation, we identify the tensor products
$k\ot V$ and $V\ot k$ with $V$. 
For a $k$-algebra $A$ we denote by $\CM_A$, 
$_A\CM$ and $_A\CM_A$ the category of  right $A$-modules, left 
$A$-modules and $A$-bimodules respectively. Similarly, for a 
$k$-coalgebra $C$  we  
 denote by $\CM^C$, $^C\!\CM$ and $^C\!\CM^C$ the category of  right 
$C$-comodules, left $C$-comodules and $C$-bicomodules respectively. 
Also, by $_A\CM^C$ ($^C\!\CM_A$) we 
denote the category of left (right) 
$A$-modules with the action ${}_V\mu$ ($\mu_V$)  and right (left) 
$C$-comodules with the coaction $\Delta_V$ (${}_V\Delta$) such that 
$\Delta_V\circ{}_V\mu = ({}_V\mu\tens C)\circ (A\tens\Delta_V)$ 
(${}_V\Delta\circ\mu_V = (C\tens\mu_V)\circ({}_V\Delta\tens A)$), 
i.e., $\Delta_V$ (${}_V\Delta$) is right (left) $A$-linear. 
For coactions and coproducts we use Sweedler's notation with suppressed 
summation sign: $\Delta_A (a) = a\sw 0\tens a\sw 1,\; \hD(c)=c\1\ot c\2\,$.

\section{\boldmath $C$-Galois extensions} 
 
First recall the definition of a Hopf-Galois extension  
(see \cite{Mon:hop} for a review). 
 
\bde\Label{hgdef} 
Let $H$ be a Hopf algebra, $A$ be a right $H$-comodule algebra, and 
$B:=A\sp{co H}:=\{a\in A\,|\; \hD_A\ a=a\ot 1\}$. We say that 
A is a (right)  {\em Hopf-Galois extension} (or $H$-Galois extension) 
of $B$ iff the canonical left 
$A$-module 
right $H$-comodule map 
$ 
can:=(m\ot H)\circ (A\ot\sb B \hD_A)\, :\; A\ot\sb B A\lra A\ot H 
$ 
is bijective. 
\ede 
 
Note that in the situation of Definition~\ref{hgdef}, both 
$A\tens_B A$ and $A\tens H$ are objects in $_A\CM^H$ via the maps 
$m\tens_B A$, $A\tens_B\Delta_A$ and $m\tens A$, $A\tens\Delta$ 
respectively. The canonical map $can$ is a morphism in this category. 
Thus the extension $B\subset A$ is Hopf-Galois if $A\tens_B A \cong 
A\tens H$ as objects in $_A\CM^H$ by the canonical map $can$. 
 
It has recently been observed in \cite{Brz:hom} that one can view 
quantum embeddable homogeneous spaces $B$ of a Hopf algebra $H$, such 
as the family (indexed by $c\in [0,\infty]$) of 
quantum two-spheres of Podle\'s   \cite{Pod:sph},  as  extensions 
by a coalgebra $C$. 
It is known (see p.200 in~\cite{Pod:sph}) that, 
 except for the North Pole sphere ($c=0$), 
no other spheres of this family are quantum quotient spaces of $SU_q(2)$. 
They escape the standard Hopf-Galois description. 
To include this important case in Galois extension theory, one needs 
to generalise the notion of a Hopf-Galois extension 
to the case of an algebra extended to an algebra by a coalgebra 
(cf.~\cite{DijKor:hom} for the dual picture). 
% rather than a Hopf algebra.  
This is obtained by weakening the requirement that $\Delta_A$ be an 
algebra map,  and leads to the notion of a coalgebra-Galois extension.  
(A special case of this kind was considered in~\cite[p.291]{Sch:nor}.)

To formulate the definition of a coalgebra\gal, first we need a general  
concept of coinvariants:
\bde[\cite{t-m97}]\Label{coinv}
Let $A$ be an algebra and a right $C$-comodule. Then
\[
A\sp{co C}:=
\left\{b\in A\; |\; \Delta_A(ba)=b\Delta_A(a),\;\forall a\in A\right\}
\]
is a subalgebra of $A$. We call it the subalgebra of 
{\em (right) coinvariants}.
\ede
Observe that when $\Delta_A$ is an algebra map, the above definition coincides
 with the usual definition of coinvariants as elements $b$ of $A$ such that
$\Delta_A(b)=b\te 1$. 
Definition~\ref{coinv} does not require the existence of
a group-like element in the coalgebra. This allows one to define 
coalgebra\gal\ for arbitrary coalgebras.
\bde\Label{cge} 
Let $C$ be a coalgebra, 
$A$ an algebra and  a right $C$-comodule and let  
$B=A\sp{co C}$. We say that 
A is a (right) {\em coalgebra-Galois extension}  
(or $C$-Galois extension) of $B$ iff the canonical left $A$-module right 
$C$-comodule map 
$
can:=(m\ot C)\ci (A\ot\sb B \hD_A)\, :
\; A\ot\sb B A\lra A\ot C\vspace*{-2.5mm} 
$ 
is bijective. 
\ede 
 
%Note that  $B$ is an algebra,  
%so that the condition of Definition~\ref{cge} makes sense. 
In what follows, we will consider only right coalgebra-Galois extensions,  
and skip writing ``right" for brevity.  
The conditions of Definition~\ref{cge} suffice to make both $A\tens_B A$ 
 and $A\tens C$ objects in $_A\CM^C$ via the maps 
$m\tens_B A$, $A\tens_B\Delta_A$ and $m\tens A$, $A\tens\Delta$,  
respectively. The canonical map is again a morphism in $_A\CM^C$. The 
extension $B\subset A$ is $C$-Galois if $can$ is an isomorphism in $_A\CM^C$. 
 
By the reasoning as in the proof of 
Proposition~1.6 in~\cite{h}, one can obtain an alternative (differential) 
definition of a \cge: 
\bpr\Label{equi1} 
Let $C$ and $B\inc A$ be as above. Let $\hO^1\! A:=\mbox{\em Ker}\, m$ denote 
the universal differential calculus on $A$, and $C^+:=\mbox{\em Ker}\,\eps$ 
the augmentation  
ideal of $C$. The algebra $A$ is 
a $C$\gal\ {\em if and only if} the following sequence of left $A$-modules 
is exact: 
\beq\Label{dif} 
0\lra A(\hO^1\! B)A\lra\hO^1\! A\st{\overline{can}}{\lra} A\ot C^+\lra 0\, , 
\eeq 
where $\overline{can}:=(m\ot C)\ci(A\ot\hD_A)$. 
\epr 
In the Hopf-Galois case, one can generalise the above sequence to  
a non-universal differential calculus in a straightforward manner: 
\beq\Label{seq} 
0\lra A\hO^1(B)A\lra\hO^1(A)\st{\tilde\chi}{\lra} A\ot (H^+/R_H)\lra 0\, , 
\eeq 
where $R_H$ is the $ad_R$-invariant right 
ideal of the Hopf algebra $H$ defining a bicovariant calculus  
on~$H$~\cite{wor}, 
and $\tilde\chi$ is defined by a formula fully analogous to the formula 
for $\overline{can}$. 
Sequence (\ref{seq}) is a starting point of the quantum-group gauge theory  
proposed in~\cite{BrzMa:gau} and continued in~\cite{h}. The Hopf-Galois  
extension describes a quantum principal bundle with the universal 
differential  
calculus. Proposition~\ref{equi1} shows that the $C$-Galois 
extension can also be viewed as a generalisation of such a bundle, 
a principal coalgebra bundle. The theory of coalgebra bundles and 
connections on them (also for non-universal differential calculus) 
was developed in \cite{BrzMa:coa}. More specifically, the theory 
considered in \cite{BrzMa:coa} uses the notion of an entwining structure 
(closely connected with the theory of factorisation of 
algebras considered in \cite{Maj:phy}) and identifies coalgebra 
principal bundles with $C$-Galois extensions 
constructed within this entwining structure.  
 
The aim of this section 
is to show that to each $C$-Galois extension of Definition~\ref{cge} there 
corresponds a natural entwining structure. Therefore the notions of a 
$C$-Galois extension of Definition~\ref{cge} and a $\psi$-principal coalgebra 
bundle of \cite[Proposition~2.2]{BrzMa:coa} are equivalent to each 
other provided that there exists a group-like $e\in C$ such that  
$\Delta_A(1) = 1\otimes e$. First we recall the definition of an  
entwining structure. 
\bde\Label{ent} 
Let $C$ be a coalgebra, $A$ an algebra and let $\psi$ be a 
$k$-linear map $\psi: C\tens A\to A\tens C$ such that 
\begin{equation}\Label{diag.A} 
\psi\circ(C\tens m) = (m\tens C)\circ (A\tens\psi)\circ(\psi\tens A), 
\qquad \psi\circ (C\tens\eta) = \eta\tens C, 
\end{equation} 
\begin{equation}\Label{diag.B} 
(A\tens\Delta)\circ\psi = (\psi\tens 
C)\circ(C\tens\psi)\circ(\Delta\tens A), \qquad (A\tens \eps)\circ\psi = 
\eps\tens A, 
\end{equation} 
where $\eta$ is the unit map $\eta: \alpha\mapsto\alpha 1$. 
Then $C$ and $A$ are said to be {\em entwined} by $\psi$ and the triple 
$(A,C,\psi)$ is called an {\em entwining structure}. 
\ede 
Entwining structures can be also understood as follows. Given an 
algebra $A$ and a coalgebra $C$ we consider $A\tens C$ as an object in 
$_A\CM$ with the structure map $m\tens C$. Similarly we consider 
$C\tens A$ as an object in $^C\!\CM$ via the 
map $\Delta\tens A$. Then we have (see also \cite[Theorem~5.2]{Sch:yet}) 
\bpr\Label{entw} 
Let $\mu_{A\tens C} : A\tens C\tens A\to A\tens C$ and 
$\Delta_{C\tens A} :C\tens A\to C\tens A\tens C$ 
be the maps making $A\tens C$ an 
object in ${}_A\CM_A$ and $C\tens A$ an object in ${}^C\!\CM^C$ 
correspondingly, and such that 
$$ 
(\eps\tens A\tens 
C)\circ\Delta_{C\tens A} = \mu_{A\tens C}\circ(\eta\tens C\tens A). 
$$ 
Then the pairs $(\mu_{A\tens C},\Delta_{C\tens A})$ of such compatible 
maps are in one-to-one correspondence with 
the entwining structures $(A,C,\psi)$.  
\note{Entwining structures $(A,C,\psi)$ are in bijective correspondence with 
the maps $\mu_{A\tens C} : A\tens C\tens A\to A\tens C$ and 
$\Delta_{C\tens A} :C\tens A\to C\tens A\tens C$ making $A\tens C$ an 
object in ${}_A\CM_A$ and $C\tens A$ an object in ${}^C\!\CM^C$ 
correspondingly, and such that 
$$ 
(\eps\tens A\tens 
C)\circ\Delta_{C\tens A} = \mu_{A\tens C}\circ(\eta\tens C\tens A). 
$$} 
\epr 
\proof First assume that $C\tens A \in {}^C\!\CM^C$ and $A\tens C \in 
{}_A\CM_A$. Then 
$$ 
(\Delta\tens A\tens C)\circ\Delta_{C\tens A} = (C\tens \Delta_{C\tens 
A})\circ (\Delta\tens A), 
$$ 
and 
$$ 
\mu_{A\tens C}\circ(m\tens C\tens A) 
= (m\tens C)\circ(A\tens \mu_{A\tens C}). 
$$ 
Define 
$$ 
\psi = (\eps\tens A\tens 
C)\circ\Delta_{C\tens A} = \mu_{A\tens C}\circ(\eta\tens C\tens A). 
$$ 
Then 
\begin{eqnarray*} 
(m\tens C)\!\!\!\!&\circ&\!\!\!\!(A\tens\psi)\circ(\psi\tens A)\\ & = & (m\tens C)\circ 
(A\tens \mu_{A\tens C})\circ(A\tens\eta\tens C\tens A)\circ(\psi\tens 
A)\\ 
& = & \mu_{A\tens C}\circ(m\tens C\tens A)\circ(A\tens\eta\tens 
C\tens A)\circ(\psi\tens 
A)\\ 
& = & \mu_{A\tens C}\circ (\mu_{A\tens C}\tens A)\circ(\eta\tens 
C\tens A\tens A)\\ 
& = & \mu_{A\tens C}\circ (A\tens C\tens m)\circ(\eta\tens 
C\tens A\tens A)\\ 
& = & \mu_{A\tens C}\circ(\eta\tens C\tens A)\circ (C\tens m) = 
\psi\circ(C\tens m). 
\end{eqnarray*} 
Furthermore 
$$ 
\psi\circ(C\tens\eta) =  \mu_{A\tens C}\circ(\eta\tens C\tens A) 
\ci(C\tens \eta) 
 = \mu_{A\tens C}\circ(\eta\tens C\tens \eta) = \eta\tens C . 
$$ 
Therefore $\psi$ satisfies conditions (\ref{diag.A}). Dualising the 
above calculation (namely, interchanging $\Delta$ with $m$, $C$ with 
$A$, $\eps$ with $\eta$ and $\Delta_{C\tens A}$ with $\mu_{A\tens 
C}$)  one easily finds that $\psi$ satisfies conditions 
(\ref{diag.B}) too. Hence $(A,C,\psi)$ is an 
entwining structure. 
 
Conversely, let $(A,C,\psi)$ be an entwining structure. Define 
$\Delta_{C\tens A} = (C\tens\psi)\circ(\Delta\tens A)$ and 
$\mu_{A\tens C} = (m\tens C)\circ (A\tens\psi)$. Then  
\begin{eqnarray*} 
(C\tens A\tens\Delta)\!\!\!\!&\circ&\!\!\!\!\Delta_{C\tens A}\\ &=& (C\tens 
A\tens\Delta)\circ(C\tens\psi)\circ(\Delta\tens A) \\ 
&=& (C\tens\psi\tens 
C)\circ(C\tens C\tens\psi)\circ(C\tens \Delta\tens A)\circ(\Delta\tens 
A)\\ 
& = & (C\tens\psi\tens 
C)\circ(C\tens C\tens\psi)\circ(\Delta\tens C\tens A)\circ(\Delta\tens 
A)\\ 
& = & (C\tens\psi\tens C)\circ(\Delta\tens A\tens 
C)\circ(C\tens\psi)\circ(\Delta\tens A)\\ 
& = &(\Delta_{C\tens A}\tens 
C)\circ\Delta_{C\tens A}. 
\end{eqnarray*} 
We used property (\ref{diag.B}) to derive the second equality, 
and then coassociativity of the coproduct to derive the third one. 
Similarly, 
\begin{eqnarray*} 
(C\tens A\tens\eps)\circ\Delta_{C\tens A}& =& (C\tens 
A\tens\eps)\circ(C\tens\psi)\circ(\Delta\tens A)\\ 
& =& (C\tens\eps\tens 
A)\circ (\Delta\tens A) = C\tens A. 
\end{eqnarray*} 
Hence $\Delta_{C\tens A}$ is a right coaction. By dualising the 
above argument one verifies that $\mu_{A\tens C}$ is a right 
action. Finally, an elementary calculation shows that 
$$ 
(\eps\tens A\tens 
C)\circ\Delta_{C\tens A} = \mu_{A\tens C}\circ(\eta\tens C\tens A) = \psi. 
$$ 
Thus the bijective correspondence is established. 
\epf 

To each entwining structure $(A,C,\psi)$ one can associate the category 
$\CM_A^C(\psi)$ of right $(A,C,\psi)$-modules, introduced and 
studied in \cite{Brz:mod}. The objects of $\CM_A^C(\psi)$ are 
right $A$-modules and right $C$-comodules $V$ such that for all $v\in V$ 
and $a\in A$, 
$$%\beq\Label{psi.module} 
\Delta_V(v\cdot a) = v\sw 0\psi(v\sw 1\tens a). 
$$%\eeq 
The morphisms in $\CM_A^C(\psi)$ are right $A$-module right $C$-comodule maps. 
Some important categories well-studied in the Hopf algebra theory, such as 
the category of right $(A,H)$-Hopf modules for a right $H$-module algebra $A$  
\cite{Doi:str},  
or the category of right-right Yetter-Drinfeld  
modules \cite{Yet:rep} can be seen as examples of $\CM_A^C(\psi)$. 
 
\note{
As a consequence of the definition of an entwining structure one easily finds 
(cf. \cite[Proposition 2.2]{BrzMa:coa}) that for any  
linear map $\bar{\eta} : k\to A\otimes C$   
such that  
$$ 
%\beq\Label{eta.1} 
(m\otimes C^{\otimes 2})\circ (A\otimes\psi\otimes  
C)\circ(\bar{\eta}\otimes\bar{\eta}) =  
(A\otimes\Delta)\circ\bar{\eta}, \quad 
(A\tens\eps)\circ\bar{\eta} =  
\eta, 
%\eeq 
$$ 
 the algebra $A$ is a right $C$-comodule    
with the coaction $\Delta_A =(m\tens   
C)\circ(A\tens\psi)\circ(\bar{\eta}\tens A)$. 
Explicitly, the coaction reads $\Delta_A(a) =   
\sum_ia^i\psi(c^i\tens a)$, where $\sum_i a^i\otimes c^i =   
\bar{\eta}(1)$. In particular $\sum_i a^i\otimes c^i = \Delta_A (1)=  
1\sw 0\otimes 1\sw 1$, so that the coaction can be equivalently written 
$\Delta_A (a) = 1\sw 0\psi(1\sw 1\otimes a)$. In this case $A$ itself is a 
right $(A,C,\psi)$-module with structure maps $m$ and $\Delta_A$. This  
can be verified by a  
direct computation which uses the definition of $\Delta_A$ and the first of  
properties (\ref{diag.A}).} 
\note{Now we can consider  
\bde\Label{psi} 
Let $(A,C,\psi)$ be an entwining structure,  $\bar{\eta}:  
k\to A\otimes C$ satisfy (\ref{eta.1}), and  $\Delta_A = (m\tens   
C)\circ(A\tens\psi)\circ(\bar{\eta}\tens A)$ be a right 
$C$-coaction on $A$. Let $B := \{b\in A \; | \; \forall a \in A, \; 
\Delta_A(ba) 
 = b\Delta_A(a)\}$. We say that $B\subset A$ is a {\em $\psi$-Galois 
 extension} of 
$B$ by a coalgebra $C$ iff $B\subset A$ is a \cge, i.e., 
iff the canonical homomorphism $can  : 
A\tens_BA\to A\tens C$ in $_A\CM^C$ is an isomorphism. 
(We call the underlying quantum-space structure of $(A,B,C,\psi,\bar{\eta})$  
 a $\psi$-principal bundle.) 
\ede 
Definition~\ref{psi} generalises the notion introduced in  
\cite[Proposition~2.2]{BrzMa:coa}, where the existence of} 
\note{In particular if there exists  a 
group-like $e\in C$ then the map $\bar{\eta}$ can defined by  
$\bar{\eta}(1) = 1\otimes e$. } 
\note{To see that with this choice of $\bar{\eta}$, 
 Definition~\ref{psi} 
defines the same object as \cite[Proposition~2.2]{BrzMa:coa}, 
it is necessary to note that in the case of the $\psi$-Galois extension,  
$B$ can be  identified with  
the set $A_0 := \{ a\in A \; | \; \Delta_A(a) = a1\sw 0\otimes 1\sw 1\}$. 
Indeed, the 
inclusion $B\subseteq A_0$ is obvious, while the fact that every $a\in A_0$ 
is an element of $B$ follows from the fact that for all $a,a'\in A$,  
$\Delta_A(aa')= a\sw 0\psi(a\sw 1\otimes a')$. 
The latter can be verified by a  
direct computation which uses the definition of $\Delta_A$ and the first of  
properties (\ref{diag.A}). Thus if $\bar{\eta}(1) = 1\otimes e$ is chosen, 
then $B = A_e^{co C}:= \{b\in A \; | \; \Delta_A(b) = b\otimes e\}$,  
as in \cite{BrzMa:coa}.
} 
 
The main result of the paper is contained in the following: 
\bth\Label{main} 
Let $A$ be a $C$-Galois extension of $B$. Then there exists a unique 
map $\psi: C\tens A\to A\tens C$ entwining $C$ with $A$ and such that  
\note{for  
all $a,a'\in A$, $\Delta_A(aa')= a\sw 0\psi(a\sw 1\tens a')$, i.e. such that}  
$A\in \CM_A^C(\psi)$ with the
structure maps $m$ and $\Delta_A$. (The map $\psi$  
is called the {\em canonical entwining map} associated to the $C$-Galois  
extension $B\subseteq A$.) 
\ethe 
\proof 
Assume that $B\subset A$ is a \cge. Then 
$can$ is bijective and there exists the 
{\it translation map} $\tau : C\to A\tens_B A$, $\tau(c):= 
can^{-1} (1\tens c)$. We use the notation $\tau(c) = c\su 
1\tens c\su 2$ (summation understood). 
Using \cite[Proposition 3.9]{Brz:tra} or 
\cite[Remark~3.4]{Sch:rep} and their obvious 
generalisation to the present case, for all $c\in C$ and $a\in A$,  
one obtains that 
 
\hspace*{2.5mm} (i) $ c\su 1c\su 2 =\eps (c)$, 
 
\hspace*{2.5mm} (ii) $ a\sw 0a\sw 1\su 1\tens a\sw 1\su 2 = 1\tens a$, 
 
\hspace*{2.5mm} (iii) $ c\su 1\tens c\su 2\sw 0\tens c\su 2\sw 1 =  c\sw 1\su 
1 \tens c\sw 1\su 2\tens c\sw 2$.

Define a map $\psi : C\tens A\to A\tens C$ by 
\beq\Label{psi1} 
\psi = can 
\circ(A\otimes_Bm)\circ(\tau \tens A),\;\;\; 
\psi (c\tens a)=  c\su 1(c\su 2 a)\sw 0\tens (c\su 2 a)\sw 1\, . 
\eeq 
We show that $\psi$ entwines $C$ and $A$. 
By definition of the translation map, 
we have: 
$$ 
\psi(c\tens 1) =  c\su 1c\su 2\sw 0\tens c\su 2\sw 1 = 1\tens c. 
$$\parindent=5mm 
Furthermore, by property (i) of the translation map, 
$$ 
\llp(\id\tens\eps)\circ\psi\lrp(c\tens a) =  c\su 1(c\su 2 a)\sw 0\tens 
\eps((c\su 2 a)\sw 1) =  c\su 1c\su 2 a = \eps(c) a. 
$$ 
Thus we have proven that the second equations  in 
(\ref{diag.A})--(\ref{diag.B}) hold. Next we have 
\begin{eqnarray*} 
\llp(m\tens C)\!\!\!\!&\ci&\!\!\!\!(A\tens\psi)\circ(\psi\tens 
A)\lrp(c\tens a\tens a')\\  
&= & 
\llp(m\tens C)\ci(A\tens\psi)\lrp 
(c\su 1(c\su 2 a)\0\tens(c\su 2 a)\1\tens a')\\ 
&= & 
c\su 1(c\su 2 a)\0(c\su 2 a)\1\su 1((c\su 2 a)\1\su 2 a')\0 
\tens((c\su 2 a)\1\su 2 a')\1\\ 
&= & 
c\su 1(c\su 2 aa')\sw 0\tens (c\su 2 aa')\sw 1 
=  
\llp\psi\circ(C\tens m)\lrp(c\tens a\tens a'), 
\end{eqnarray*} 
where we used the property (ii) of the translation map to derive the 
third equality. 
Hence  the first of equations (\ref{diag.A}) holds. 
Similarly, 
\begin{eqnarray*} 
\llp(\psi\tens C)\!\!\!\!&\ci&\!\!\!\!(C\tens\psi)\circ(\Delta\tens A)\lrp (c\tens a) \\ 
&=& 
(\psi\tens C)( c\sw 1 \tens c\sw 2\su 1(c\sw 2\su 2 
a)\sw 0\tens (c\su 2 a)\sw 1) \\ 
&=& 
  c\sw 1\su 1(c\sw 1\su 2 c\sw 2\su 1(c\sw 2\su 2a)\sw 0)\sw 0 
\tens (c\sw 1\su 2 c\sw 2\su 1(c\sw 2\su 2a)\sw 0)\sw 1 \\ 
&&\tens (c\sw 
2\su 2a)\sw 1\\ 
&=& 
 c\su 1(c\su 2\sw 0 c\su 2\sw 1\su 1(c\su 2\sw 1\su 2a)\sw 0)\sw 0 
\tens (c\su 2\sw 0 c\su 2\sw 1\su 1(c\su 2\sw 1\su 2a)\sw 0)\sw 1 
\\ 
&&\tens (c\su 2\sw 1\su 2 a)\sw 1\\ 
&= &c\su 1((c\su 2 
a)\sw 0)\sw 0\tens ((c\su 2 a)\sw 0)\sw 1\tens (c\su 2 a)\sw 1 \\ 
&=& 
 c\su 1(c\su 2 
a)\sw 0\tens (c\su 2 a)\sw 1\tens (c\su 2 a)\sw 2  
= 
\llp(A\tens\Delta)\circ\psi\lrp(c\tens a) . 
\end{eqnarray*} 
We used property (iii) of the translation map to derive the third 
equality and then property (ii) to derive the fourth one. 
Hence $C$ and $A$ are entwined by $\psi$ as required. 

Now, using (ii),  
we have for all $a,a'\in A$ 
\begin{eqnarray*} 
a\sw 0\psi(a\sw 1\tens a') & = & 
a\sw 0a\sw 1\su 1(a\sw 1\su 2 a')\sw 0 \otimes 
(a\sw 1\su 2 a')\sw 1\\ 
& = &  (aa')\sw 0\otimes (aa')\sw 1 = \Delta_A(aa'), 
\end{eqnarray*} 
{\em i.e.} $A$ is an $(A,C,\psi)$-module with structure maps $m$ and  
$\Delta_A$. 
It remains to prove the uniqueness of the entwining map $\psi$ given by  
(\ref{psi1}). Suppose that there is an entwining map $\tilde{\psi}$ such 
that $A\in \CM_A^C(\tilde{\psi})$ with structure maps $m$ and $\Delta_A$.  
Then, for all $a\in A$, $c\in C$, 
$$ 
\psi(c\tens a) = c\su 1(c\su 2 a)\sw 0\tens (c\su 2 a)\sw 1  =  
c\su 1c\su 2\sw 0 \tilde{\psi}(c\su 2\sw 1\tens a)  
= \tilde{\psi}(c\tens a), 
$$ 
where we used the definition of the translation map to obtain the 
last equality. 
%Thus the uniqueness of $\psi$ is proven and the proof 
%of the theorem is completed. 
\epf

\section{\boldmath $A$-Galois coextensions} 
 
The dual version of a Hopf-Galois extension can be viewed as a  
non-commutative generalisation of the theory of quotients 
of formal schemes under free actions of formal group schemes  
(cf.~\cite{Sch:pri}). In this section we dualise $C$-Galois extensions 
and derive results analogous to the results discussed in the previous 
section. 

First recall the definition of a cotensor product.  
Let $B$ be a coalgebra and $M$, $N$ a right and left 
$B$-comodule respectively. 
\note{Let $C$, $B$ be 
coalgebras and $\pi:C\to B$ be a coalgebra surjection. } 
The cotensor product $M\Box_{B}N$ is defined by the exact sequence 
\[ 
0\rightarrow M\Box_{B}N\hookrightarrow M\ot N\st{\ell} 
{\rightarrow}M\ot B\ot N, 
\] 
where $\ell$ is the coaction equalising map  
$\ell=\Delta\!_M\ot N-M\ot\, _N\!\Delta$. 
\note{$$ 
C\Box_{B}C \subset C\tens C 
\stackrel{\textstyle\longrightarrow}{\longrightarrow} 
C\tens B\tens C. 
$$ 
The maps on the right hand side of the equaliser  are 
 $(C\tens\pi\tens C)\circ(C\tens\Delta)$ and 
$(C\tens\pi\tens C)\circ(\Delta\tens C)$.} 
In particular, if $C$ is a coalgebra, $I$ its coideal and $B=C/I$, then 
\[ 
C\Box_{B}C\!=\!\left\{\!\sum_i\!c^i\ot \tc^i\!\in\! C\ot C\;\left| 
\;\sum_i\! c^i\1\ot\pi(c^i\2)\ot\tc^i 
\!=\!\sum_i\!c^i\ot\pi(\tc^i\1)\ot \tc^i\2\right.\right\}, 
\] 
 where $\pi:C\to B=C/I$ is the canonical 
surjection. The following definition \cite[p.3346]{sch} 
dualises the concept of a \hge: 
\begin{definition}\Label{hgc} 
Let $H$ be a Hopf algebra, $C$ a right $H$-module coalgebra with the 
action $\mu_C:C\otimes H\to C$. Then   
$I := \{\mu_C(c,h)-\epsilon(h)c\;\;|\;\; c\in C, h\in H\} $ is 
a coideal in $C$ and thus $B:=C/I$ is a coalgebra. We say that 
$C\twoheadrightarrow B$ is a (right) {\em Hopf-Galois coextension}  
(or $H$-Galois coextension) iff the  
canonical left $C$-comodule right $H$-module map 
$ 
cocan := (C\ot\mu_C)\circ(\Delta\tens H):C\te H \ra  
C\Box\sb BC 
$ 
is a bijection. 
\end{definition} 
With the help of the property  
$\hD\ci\mu_C=(\mu_C\ot\mu_C)\ci(C\ot\mbox{flip}\ot H)\ci(\hD\ot\hD)$, 
it can be directly checked that the image of the map $cocan$ 
is indeed contained in $C\Box_BC$. 
To see more clearly that Definition~\ref{hgc} is obtained by dualising 
 Definition~\ref{hgdef}, one can notice that both 
$C\tens H$ and $C\Box\sb B C$ are objects in $\sp C\!\CM\sb H$, which 
is dual to $\sb A\CM\sp H$. The structure maps are 
$\Delta\tens C$, $C\tens m$ and $\Delta\Box\sb BC$, $C\Box\sb B\mu_C$ 
respectively. The canonical map $cocan$ is a morphism in $\sp C\!\CM\sb H$.  
The coextension $C\twoheadrightarrow B$ is Hopf-Galois if $C\tens H \cong 
C\Box_BC$ as objects in $\sp C\!\CM\sb H$ by the canonical map $cocan$. 
 
The notion of a Hopf-Galois coextension can be generalised by replacing 
$H$ by an algebra $A$ and weakening the condition that the action $\mu_C$ 
is a coalgebra map. This generalisation dualises the construction of 
a $C$-Galois extension of the previous section. First we prove 
\note{ 
\ble\Label{colem} 
 Let $C$ be a coalgebra, $A$ an algebra and $\kappa :A\to 
k$ an algebra character. Assume that $C\in \CM_A$ with the action 
$\mu_C$ and define 
\begin{equation}\Label{ikappa} 
I_\kappa = \{\mu_C(c,a)-\kappa(a)c\;\;|\;\;  c\in C, a\in A\}. 
\end{equation} 
Let 
$B:= C/I_\kappa$. If $\eps\circ\mu_C = \eps\tens\kappa$ and  
$ 
\Delta_\kappa\circ\mu_C = 
(B\otimes\mu_C)\circ(\Delta_\kappa\tens A), 
$ 
where $\Delta_\kappa = (\pi_\kappa\tens C)\circ\Delta$  
and $\pi_\kappa : C\to C/I_\kappa$ is the  
canonical surjection, then $B$ is a coalgebra. 
(If the conditions of the lemma are satisfied we say that 
the action $\mu_C$ is {\em left $B$-colinear} and 
{\em counital with respect to 
$\kappa$}.) 
\ele 
\proof By the counitality of $\mu_C$, we have $I_\kappa 
\subset C^+:=\mbox{Ker}\,\eps$. Furthermore, using the colinearity of $\mu_C$ one 
finds that $\Delta_\kappa(I_\kappa)\subset B\tens I_\kappa$. Since 
$\mbox{Ker}\,\pi_\kappa = I_\kappa\,$, we obtain $\Delta(I_\kappa) 
\subset C\tens I_\kappa + I_\kappa\tens C$. Therefore $I_\kappa$ is a 
coideal and $B$ is a coalgebra, as stated. 
\epf\ \\ 
The above lemma implies that in this case $C\in {}^B\!\CM_A$ with the 
structure maps $\mu_C$ and~$\Delta_\kappa\,$. 
}
\ble\Label{coa}
Let $A$ be an algebra and $C$ a coalgebra and right $A$-module with an action
$\mu_C:C\te A\ra C$. Then the space 
\[ 
{\rm span}\{\mu_C(c,a)\1\ha(\mu_C(c,a)\2)-c\1\ha(\mu_C(c\2,a))\;|
\; a\in A,\, c\in C,\,\ha\in\mbox{\em Hom}(C,k)\} 
\]
is a coideal of $C$. 
\ele
\bpf
Let $I$ denote the space defined in Lemma~\ref{coa}. We prove the lemma
 by showing that $D:=\{\hb\in\mbox{Hom}(C,k)\;|\; \hb(I)=0\}$
is a subalgebra of the convolution algebra $\mbox{Hom}(C,k)$ 
(see~\cite[Proposition~1.4.6 c)]{swe}). 
Note first that
$\he\in D$. Furthermore,
\bea
D=\{\hb\in\mbox{Hom}(C,k)\!\!\! &|&\!\!\! (\hb*\ha)(\mu_C(c,a))=
\hb(c\1)\ha(\mu_C(c\2,a)),
\nonumber\\
&& \fa\, a\in A,\, c\in C,\,\ha\in\mbox{Hom}(C,k)\}.\nonumber
\eea
If $\hb_1,\,\hb_2\in D$, then for any $a,\, c,\, \ha$ we have
\bea
((\hb_1*\hb_2)*\ha)(\mu_C(c,a))
& = &(\hb_1*(\hb_2*\ha))(\mu_C(c,a))
\nonumber\\ 
& = &\hb_1(c\1)(\hb_2*\ha)(\mu_C(c\2,a))
\nonumber\\ 
& = &\hb_1(c\1)\hb_2(c\2)\ha(\mu_C(c\3,a))
\nonumber\\ 
& = &(\hb_1*\hb_2)(c\1)\ha(\mu_C(c\2,a)).
\nonumber
\eea
Hence $\hb_1*\hb_2\in D$, and $D$ is a subalgebra of $\mbox{Hom}(C,k)$,
as needed.
\epf
\ble\Label{im}
Let $C$ and $A$ be as above and let $I$ denote the coideal of Lemma~\ref{coa}.
Put $B:=C/I$, and denote by $\pi:C\ra C/I$ the canonical surjection. Then the
action $\mu_C$ is left $B$-colinear, i.e., 
$(\pi\te C)\ci\hD\ci\mu_C = (B\te\mu_C)\ci((\pi\te C)\cc\hD\ot A)$,
and
$((C\te\mu_C)\cc(\hD\ot A))(C\ot A)\inc C\Box_BC$.
\ele
\bpf
Choose $a\in A,\, c\in C$. The above $B$-colinearity condition can be 
written as
$\pi(\mu_C(c,a)\1)\ot\mu_C(c,a)\2 = \pi(c\1)\ot\mu_C(c\2,a)$
It is equivalent to the condition
\[
\pi(\mu_C(c,a)\1)\ot\ha(\mu_C(c,a)\2) = \pi(c\1)\ot\ha(\mu_C(c\2,a)),\;
\fa\,\ha\in\mbox{Hom}(C,k),
\]
which can be written as
\[
\pi(\mu_C(c,a)\1\ha(\mu_C(c,a)\2) - c\1\ha(\mu_C(c\2,a))) =0,\;
\fa\,\ha\in\mbox{Hom}(C,k).
\]
This proves our first assertion.

As for the second assertion, observe that it can be stated as
\[
c\1\ot\pi(c\2)\ot\mu_C(c\3,a) =  c\1\ot\pi(\mu_C(c\2,a)\1)\ot\mu_C(c\2,a)\2\, .
\]
This is true by the left $B$-colinearity argument applied to the last two 
tensorands.
\epf

Now we can conclude that we have a well-defined map
%\beq\Label{coc} 
$$
cocan := (C\ot\mu_C)\circ(\Delta\tens A)\; :\;\; C\tens A \lra  
C\Box\sb BC, 
$$
%\eeq 
and can consider: 
\begin{definition}\Label{agc} 
Let $A$ be an algebra,  
$C$ a coalgebra and right $A$-module,  and $B=C/I$, 
where $I$ is the coideal of Lemma~\ref{coa}. We say 
that $C$ is a (right)
 {\em algebra-Galois coextension} (or $A$-Galois coextension) 
of $B$ iff the canonical left $C$-comodule right $A$-module  
map 
$
cocan := (C\ot\mu_C)\circ(\Delta\tens A):C\tens A \lra  
C\Box\sb BC
$ 
is bijective. 
\ede 
Again, to see that Definition~\ref{agc} dualises the notion of 
a $C$-Galois extension one can notice that both 
$C\tens A$ and $C\Box\sb BC$ are objects in $\sp C\!\CM\sb A$, which 
is dual to $\sb A\CM\sp C$. The structure maps are 
$\Delta\tens C$, $C\tens m$ and $\Delta\Box\sb BC$, $C\Box\sb B\mu_C$, 
respectively. The canonical map $cocan$ is a morphism in $\sp C\!\CM\sb A$.  
The right coextension $C\twoheadrightarrow B$ is $A$-Galois if $C\tens A \cong 
C\Box_BC$ as objects in $\sp C\!\CM\sb A$ by the canonical map $cocan$.
(In what follows, we consider only right coextensions and omit ``right" for
brevity.)
 
Very much as in the previous section, it turns out that every $A$-Galois 
extension 
is equipped with an entwining structure. More precisely, we have the following 
dual version of Theorem~\ref{main}: 
 
\bth\Label{main2}
Let $C$ be an $A$-Galois coextension of $B$. Then there exists 
a unique map $\psi :C\otimes A\to A\otimes C$ entwining $C$ with $A$  
and such that $C\in{\CM}\sp{C}_A(\psi)$ with the structure maps \hD\ and
$\mu_C$. (The map $\psi$  
is called the {\em canonical entwining map} associated to the $A$-Galois  
coextension $C\twoheadrightarrow B$.) 
\ethe 
\proof We dualise the proof of Theorem~\ref{main}. 
Assume that $C\twoheadrightarrow B$ is an $A$-Galois coextension. Then 
$cocan$ is a bijection and there exists the  
{\it cotranslation map} $\check{\tau} : C\Box_B C \to 
A$, $\check{\tau} := (\eps\tens A)\circ cocan^{-1}$. By dualising 
properties of the translation map  
(or directly from the definition of $\check{\tau}$), 
one can establish the following properties of the cotranslation map: 
 
\hspace*{2.5mm}(i) $\check{\tau}\ci\hD = \eta\ci\he$, 
 
\hspace*{2.5mm}(ii) $\mu_C\circ(C\tens\check{\tau})\circ(\Delta \ot C) = 
\eps\tens C$ on $C\Box_BC$, 
%or explicitly $\mu_C(c\sw 1,\check{\tau}(c\sw 2,c') = \eps(c)c'$, 
 
\hspace*{2.5mm}(iii) $\check{\tau}\ci(C\ot\mu_C)=m\ci(\check{\tau}\ot A)$
on $C\Box_BC\ot A$,
 
%$\check{\tau}(c,\mu_C(c',a)) = \check{\tau}(c,c')a$, 
\note{ 
for any $c,c'\in C$ and $a\in A$. Combining (iii) with (ii)  one obtains 
 \hspace*{7.5mm}(iv) $\check{\tau}(c,c'\sw 1)\check{\tau}(c'\sw 2,c'') =  
\eps(c')\check{\tau}(c,c'')$, for any $c,c',c''\in C$. 
} 
It follows from (iii) that
\[
m\ci(\check{\tau}\te\check{\tau})=\check{\tau}\ci(C\ot\mu_C\cc(C\te\check{\tau}))
\mbox{\ on\ }C\Box_BC\ot C\Box_BC.
\]
Consequently, we can conclude from (ii) that
\beq\Label{iv}
m\ci(\check{\tau}\te\check{\tau})\ci(C\te\hD\te C)
=\check{\tau}\ci(C\te\he\te C) \mbox{\ on\ } C\Box_BC\Box_BC.
\eeq
(Observe that $(C\te\hD\te C)(C\Box_BC\Box_BC)\inc C\Box_BC\ot C\Box_BC$.
Here and below one has to pay attention to the domains of certain mappings.)

Using the cotranslation map and noticing that $(C\ot\hD)(C\Box_BC)\inc 
C\Box_BC\ot C$,
 one defines a map $\psi :C\tens A\to A\tens C$ by 
\begin{eqnarray*} 
&&\psi = 
(\check{\tau}\tens C)\circ(C\tens\Delta)\circ cocan,\\ 
&&\psi(c\tens a ) =  \check{\tau}(c\sw 1,\mu_C(c\sw 2,a)\sw 
1)\tens\mu_C(c\sw 2,a)\sw 2\, . 
\end{eqnarray*}%\parindent=5pt 
 We now show that $\psi$ entwines 
$C$ and $A$. For any $c\in C$ we find 
$$ 
\psi(c\tens 1) = \check{\tau}(c\sw 1,\mu_C(c\sw 2,1)\sw 1)\tens \mu_C(c\sw 
2,1)\sw 2 = \check{\tau}(c\sw 1,c\sw 2)\tens c\sw 3 = 1\tens c, 
$$ 
where we used property (i) to derive the last equality. Furthermore, 
\begin{eqnarray*}
(A\te\he)\ci\psi &=& ((\he\te A)\ci cocan\sp{-1}\ot\he)\ci(C\ot\hD)\ci cocan\\
&=& (\he\te A)\ci cocan\sp{-1}\ci cocan=\he\te A.
\end{eqnarray*}
Thus we have proven that the second conditions of 
(\ref{diag.A}) and (\ref{diag.B}) are fulfilled by $\psi$. 
\note{
for any $c\in C$, $a\in A$, 
\begin{eqnarray*} 
\llp(A\te\eps)\circ\psi\lrp(c\te a)   
&=&  
\check{\tau}(c\1,\mu_C(c\sw 2,a)\1)\eps(\mu_C(c\2,a)\2)  
\\ &=& 
 \check{\tau}(c\1,\mu_C(c\2,a)) 
 = 
  \check{\tau}(c\1,c\2)a = 
 \eps(c) a, 
\end{eqnarray*} 
where the penultimate equality is implied by (iii), and the ultimate 
one by~(i). 
}
To prove the first of equations (\ref{diag.A}), 
we compute 
\begin{eqnarray} 
&&
(m\tens C)\ci(A\tens\psi)\circ(\psi\tens A)
\nonumber\\  && =
(m\cc(A\te\check{\tau})\ot C)\ci(A\ot C\ot\hD)%\ci(A\ot cocan)
\ci(\check{\tau}\ot cocan)\ci(C\ot\hD\ot A)\ci(cocan\ot A)
\nonumber\\  && =
(m\cc(\check{\tau}\te\check{\tau})\ot C)\ci(C\sp{\ot3}\ot\hD)
\ci(C\sp{\ot2}\ot cocan)\ci(C\ot\hD\ot A)\ci(cocan\ot A)
\nonumber\\  && =
(m\cc(\check{\tau}\te\check{\tau})\ot C)\ci(C\sp{\ot3}\ot\hD)
\ci(C\ot (C\te cocan)\cc(\hD\te A))\ci(cocan\ot A)
\nonumber\\  && =
(m\cc(\check{\tau}\te\check{\tau})\ot C)
\ci
(C\sp{\ot3}\ot\hD)
\ci
(C\ot(\hD\te C)\cc cocan)
\ci
(cocan\ot A)
\nonumber\\  && =
(m\cc(\check{\tau}\te\check{\tau})\cc(C\te\hD\te C)\ot C)
\ci
(C\sp{\ot2}\ot\hD)
\ci
(C\ot cocan)
\ci
(cocan\ot A).\nonumber
\eea
Taking advantage of (\ref{iv}), we obtain
\bea\Label{v}
&&
(m\tens C)\ci(A\tens\psi)\circ(\psi\tens A)
\nonumber\\  && =
(\check{\tau}\cc(C\te\he\te C)\ot C)
\ci
(C\sp{\ot2}\ot\hD)
\ci
(C\ot cocan)
\ci
(cocan\ot A).
\nonumber\\  && =
(\check{\tau}\ot C)
\ci
\llp C\ot (\he\te C\sp{\ot2})\cc(C\te\hD)\cc cocan\lrp
\ci
(cocan\ot A).
\eea
On the other hand, for any $c\in C,\, a\in A$, we have
\bea\Label{vi}
\llp(\he\ot C\sp{\ot2})\ci(C\ot\hD)\ci cocan\lrp(c\ot a)
\nonumber & = &
\mu_C(c,a)\1\ot\mu_C(c,a)\2
\nonumber\\ & = &
(\hD\ci\mu_C)(c\ot a).
\eea
Combining this with (\ref{v}) yields
\bea\Label{vii}
(m\tens C)\ci(A\tens\psi)\circ(\psi\tens A)
\nonumber & = &
(\check{\tau}\ot C)
\ci
(C\ot\hD)
\ci
(C\ot\mu_C)
\ci
(cocan\ot A)
\nonumber\\  & = &
(\check{\tau}\ot C)
\ci
(C\ot\hD)
\ci
cocan
\ci
(C\ot m)
\nonumber\\  & = &
\psi\ci(C\ot m), \nonumber
\note{ 
\llp(m\tens C)\circ(A\tens\psi)\lrp 
\llp\check{\tau}(c\1,\mu_C(c\sw 2,a)\1)\tens\mu_C(c\2,a)\2 \tens a'\lrp\\ 
&=&  
\check{\tau}\llp c\1,\mu_C(c\2,a)\1\lrp 
\check{\tau}\llp\mu_C(c\2,a)\2,\mu_C(\mu_C(c\2,a)\3,a')\1\lrp\\ 
&&\tens\mu_C\llp\mu_C(c\sw 2,a)\sw 3,a'\lrp\2\\ 
&=& 
\check{\tau}\llp c\1,\mu_C(\mu_C(c\2,a),a')\1\lrp 
\tens\mu_C\llp\mu_C(c\2,a),a'\lrp\2\\ 
&=& 
\check{\tau}\llp c\1,\mu_C(c\2,aa')\1\lrp\tens\mu_C(c\2,aa')\2  
=  
\psi(c\tens aa').
} 
\end{eqnarray} 
as desired. Here we used the property that $\mu_C$ is an action to derive
the penultimate equality.
\note{
We used property (iv) of the cotranslation map to derive the third 
equality. The fourth equality follows from the fact that $\mu_C$ is a 
right action. 
}
To prove the first of equations (\ref{diag.B}), first we observe that 
\beq\Label{ccc}
(\hD\ot C)\ci(cocan)=(C\ot cocan)\ci(\hD\ot A).
\eeq
Hence we obtain
% we take any $c\in C$, $a\in A$, use (ii) and compute 
\begin{eqnarray*} 
&&
(\psi\tens C)\circ(C\tens\psi)\ci(\Delta\tens A)\\
&& = 
(\psi\te C)
\cc
(C\te\check{\tau}\te C)
\cc
(C\sp{\ot2}\te\hD)
\cc
(C\te cocan)
\cc
(\hD\te A)
\\ && = 
(\psi\te C)
\cc
(C\te\check{\tau}\te C)
\cc
(C\sp{\ot2}\te\hD)
\cc
(\hD\te C)
\cc
(cocan)
\\ && = 
(\psi\te C)
\cc
(C\te\check{\tau}\te C)
\cc
(\hD\te C\sp{\ot2})
\cc
(C\te\hD)
\cc
(cocan)
\\ && = 
(\check{\tau}\te C\sp{\ot2})
\cc
(C\te\hD\te C)
\cc
\llp cocan\cc(C\te\check{\tau})\cc(\hD\te C)\te C\lrp
\cc
(C\te\hD)
\cc
(cocan).
\note{
&= & 
\llp\psi\tens C\lrp 
\llp c\1\tens\check{\tau}(c\2,\mu_C(c\3,a)\1)\tens\mu_C(c\3,a)\2\lrp\\ 
&= & 
\check{\tau}\llp c\1,\mu_C(c\2,\check{\tau}(c\3,\mu_C(c\4,a)\1))\1\lrp\\ 
&&\tens\mu_C\llp c\2,\check{\tau}(c\3,\mu_C(c\4,a)\1)\lrp\2  
\tens\mu_C(c\4,a)\2\\ 
&= & 
\check{\tau}\llp c\1,\mu_C(c\2,a)\1\lrp 
\tens\mu_C(c\2,a)\2 
\tens\mu_C(c\2,a)\3\\ 
&=& 
\llp(A\tens\Delta)\circ\psi\lrp(c\tens a).
} 
\end{eqnarray*}  
To finish the calculation, we note that 
$(C\ot\check{\tau})\ci(\hD\ot C)=cocan\sp{-1}$. Indeed, thanks to 
(\ref{ccc}), we have
\begin{eqnarray*}
(C\ot\check{\tau})\ci(\hD\ot C)\ci cocan
& = &
(C\ot\check{\tau})\ci(C\ot cocan)\ci(\hD\ot A)
\\ & = &
(C\ot\he\ot A)\ci(\hD\ot A)
\\ & = & C\ot A.
\end{eqnarray*}
\note{
Next observe  that, by 
the counitality of $\mu_C$, we have $(\eps\tens\eps)\circ cocan = 
(\eps\tens\kappa)$, whence $\kappa\circ\check{\tau} = 
\eps\tens\eps$. This applied to the definition of $\psi$ implies that 
$(\kappa\tens C)\circ\psi = \mu_C$, as needed. 
Finally, the 
uniqueness of $\psi$ can be verified by the dualisation of the  argument at 
the end of the proof of Theorem~\ref{main}. (To perform 
explicit calculations one needs to use the first of equations 
(\ref{diag.B}) and the property (iii) of the cotranslation map.) 
}
Hence
\begin{eqnarray*}
(\psi\tens C)\circ(C\tens\psi)\ci(\Delta\tens A)
& = &
(\check{\tau}\te C\sp{\ot2})
\cc
(C\te\hD\te C)
\cc
(C\te\hD)
\cc
(cocan)
\\ & =&
(\check{\tau}\te C\sp{\ot2})
\cc
(C\sp{\ot2}\te\hD)
\cc
(C\te\hD)
\cc
(cocan)
\\ & =&
(A\te\hD)
\cc
(\check{\tau}\te C)
\cc
(C\te\hD)
\cc
(cocan)
\\ & = & (A\te\hD)\cc\psi,
\end{eqnarray*}
as needed. Thus we have proved that $(A,C,\psi)$ is an entwining structure.

The next step is to show that $C\in{\CM}_A\sp{C}(\psi)$, i.e., 
$\hD\ci\mu_C=(\mu_C\ot C)\ci(C\ot\psi)\ci(\hD\ot A)$. With the help
of (\ref{ccc}), property (ii) of the cotranslation map $\check{\tau}$ 
and then (\ref{vi}), we compute:
\begin{eqnarray*}
&&
(\mu_C\ot C)
\ci
(C\ot\psi)
\ci
(\hD\ot A)
\\ && =
(\mu_C\ot C)
\ci
(C\ot\check{\tau}\ot C)
\ci
(C\sp{\ot2}\ot\hD)
\ci
(C\ot cocan)
\ci
(\hD\ot A)
\\ && =
\llp\mu_C\cc(C\te\check{\tau})\ot C\lrp
\ci
(C\sp{\ot2}\ot\hD)
\ci
(\hD\ot C)
\ci
cocan
\\ && =
\llp\mu_C\cc(C\te\check{\tau})\cc(\hD\te C)\ot C\lrp
\ci
(C\ot\hD)
\ci
cocan
\\ && =
(\he\ot C\sp{\ot2})
\ci
(C\ot\hD)
\ci
cocan
\\ && = \hD\ci\mu_C\, .
\end{eqnarray*}
As for the uniqueness of $\psi$, suppose that there exists another entwining
map $\widetilde{\psi}$ such that $C\in{\CM}_A\sp{C}(\widetilde{\psi})$.
Then
\begin{eqnarray*}
\psi & = &
(\check{\tau}\ot C)
\ci
(C\ot\hD)
\ci
(C\ot\mu_C)
\ci
(\hD\ot A)
\\ & = &
(\check{\tau}\ot C)
\ci
\llp C\ot(\mu_C\te C)\cc(C\te\widetilde{\psi})\cc(\hD\te A)\lrp
\ci
(\hD\ot A)
\\ &=&
(\check{\tau}\cc(C\te\mu_C)\ot C)
\ci
(C\sp{\ot2}\ot\widetilde{\psi})
\ci
(C\ot\hD\ot A)
\ci
(\hD\ot A)
\\ &=&
(\check{\tau}\cc(C\te\mu_C)\ot C)
\ci
(\hD\ot\widetilde{\psi})
\ci
(\hD\ot A)
\\ & =&
(\check{\tau}\cc(C\te\mu_C)\cc(\hD\te A)\ot C)
\ci
(C\ot\widetilde{\psi})
\ci
(\hD\ot A)
\\ & =&
(\check{\tau}\cc cocan\ot C)
\ci
(C\ot\widetilde{\psi})
\ci
(\hD\ot A)
\\ & =&
(\he\ot A\ot C)
\ci
(C\ot\widetilde{\psi})
\ci
(\hD\ot A) = \widetilde{\psi}.
\end{eqnarray*}
\epf

\section{Galois entwining structures}

Theorem~\ref{main} allows one to view a $\psi$-principal bundle 
as a \cge.  
More precisely, recall from \cite{BrzMa:coa} the following 
\bde\Label{psi} 
Let $(A,C,\psi)$ be an entwining structure, and let $e\in C$ be a group-like  
element. Then $B := \{b\in A \; | \; \psi(e\otimes b) = b\otimes e\}$ is an 
algebra, and we say that  
$A(B,C,\psi,e)$ is a {\em coalgebra $\psi$-principal  bundle} iff the map $can_\psi 
: A\otimes _B A\to A\otimes C$, $a\otimes a'\mapsto a\psi(e\otimes a')$ 
is bijective. 
\ede 
\bpr\Label{cor} 
For a given entwining structure $(A,C,\psi)$ and a group-like $e\in C$, the 
 following statements are equivalent:\medskip \\  
\hspace*{2.5mm} 
(1) $A(B,C,\psi,e)$ is a coalgebra $\psi$-principal bundle.\medskip \\  
\hspace*{2.5mm} 
(2) There exists a unique coaction $\hD_A:A\ra A\ot C$ such that 
$B\subset A$ is a \cge\ of $B$ by $C$, $\psi$ is the canonical entwining 
map,  and $\Delta_A(1) = 1\otimes e$. 
\epr 
\proof (1) $\Rightarrow$ (2). By \cite[Proposition~2.2]{BrzMa:coa}, $A$ is 
a right $C$-comodule with the coaction $\Delta_A : a\mapsto \psi(e\otimes a)$. 
Using the second of conditions (\ref{diag.A}) one verifies that 
$\Delta_A(1) = 1\otimes e$. The first of conditions (\ref{diag.A}) implies 
that 
$A\in \CM_A^C(\psi)$ with the structure maps $\Delta_A$ and $m$. Therefore, 
by \cite[Lemma~3.3]{Brz:mod}, $B$ as defined in Definition~\ref{psi} is 
identical with the set 
$\{b\in A\; |\; \forall a\in A, \; \hD_A(ba)=b\hD_A(a)\}$ 
considered in Definition~\ref{coinv}. Hence $\can = \can_\psi$ and 
$B\subseteq A$ is a $C$-Galois extension. By Theorem~\ref{main},
$\psi$ is the associated canonical entwining map. 
 Assume now that there exists another coaction $\hD_A'$ 
satisfying condition~(2).
Then, by Theorem~\ref{main}, $A\in\CM_A^C(\psi)$ with the structure maps 
$\hD_A'$ and $m$. Consequently, as $\hD_A'(1)=1\ot e$,
we have 
\[
\hD_A'(a)=\hD_A'(1\cdot a)=\llp(m\ot C)\ci(A\ot\psi)\lrp(\hD_A'(1)\ot a)
=\psi(e\ot a)=\hD(a),\;\;\;\fa\; a\in A.
\]

(2) $\Rightarrow$ (1). Since $\psi$ is the canonical entwining map, $A\in 
\CM_A^C(\psi)$ by Theorem~\ref{main},
so that \cite[Lemma~3.3]{Brz:mod} implies that $B$ as 
defined in Definition~\ref{coinv} is identical with the set $\{b\in A \; |\;  
\psi(e\otimes b) = b\otimes e\}$ as required in Definition~\ref{psi}.  
Furthermore, since the normalisation condition 
$\Delta_A(1) = 1\otimes e$ implies that $\tau(e) = 1\otimes _B 1$, we have
$\psi(e\otimes a) = e\su 1(e\su 2 a)\sw 0\ot (e\su 2 a)\sw 1 = \Delta_A(a)$,
for all $a\in A$. Therefore 
$can_\psi = can$ and $can_\psi$ is bijective as required. 
\endproof 
\bre\em 
\rm A slightly different definition of a $C$-Galois extension was proposed 
in \cite{Brz:hom}. Let $A$ be an algebra, $C$ a coalgebra 
and $e$ a group-like element of $C$. 
One assumes that  $A\tens C \in \CM_A\,$,  $A\in 
{}\CM^C$, and the action and coaction are such that $\Delta_A\circ 
m = \mu_{A\tens C}\circ(\Delta_A\tens A)$ and $\mu_{A\tens C}(a\tens 
e,a') = aa'\sw 0\tens a'\sw 1$ for any $a,a'\in A$. Then 
$B:= \{b\in A \; | \; \Delta_A(b) = b\otimes e\}$ is an algebra and   
the canonical map $can :A\tens_BA\to A\tens C$ 
is well-defined. $B\subset A$ is a $C$-Galois extension if the canonical 
map $can$ is a bijection. One easily finds, however, that given $A$ and  
$C$ satisfying the above conditions, $B = \{ b\in A \; | \; \forall a\in 
A, \; \Delta_A(ba) = b\Delta_A(a)\}$, so that $A\in {}_B\CM^C$ via the 
maps $\mu_A 
= m$ and $\Delta_A$. Hence, by Theorem~\ref{main}, this 
definition of a $C$-Galois extension is equivalent to the one introduced 
in~\cite{BrzMa:coa} and in Definition~\ref{cge} provided that $\Delta_A(1) = 
1\otimes e$. 
\ere 
\bre\Label{psirem}\em 
In the Hopf-Galois case, the formula for $\psi$ becomes quite simple: 
$\mbox{$\psi(h\te a)$}=a\0\ot ha\1\,$.  
If the Hopf algebra $H$ has a bijective 
antipode, $\psi$ is an isomorphism,  
and its inverse is given by 
$\psi^{-1}(a\ot h):=hS^{-1}(a\1)\ot a\0\,$.  
(In fact, $\psi$ is an isomorphism 
{\em if and only if} $H$ has a bijective antipode 
\cite[Theorem~6.5]{Brz:mod}.)
Furthermore, the coaction 
$_A\hD:A\to \mbox{$H^{op}\ot A$}$,  
$$ 
_A\hD(a):=a_{(-1)}\ot a\0:=S^{-1}(a\1)\ot a\0 
$$ 
makes $A$ a left $H^{op}$-comodule algebra, 
where $H^{op}$ stands for the Hopf algebra with the opposite multiplication. 
  One can define the following 
left version of the canonical map: $can_L:A\ot_BA\rightarrow H\sp{op}\ot A$,  
$can_L(a\ot_Ba'):=a_{(-1)}\ot a\0 a'$. It is straightforward to verify that 
$\psi\ci can_L=can$. 
Since $\psi$ and $can$ are isomorphisms, we can immediately 
conclude that so is $can_L\,$.  
\ere 
 
In parallel to the theory of coalgebra
$\psi$-principal bundles, the notion of a dual $\psi$-principal bundle was
introduced in~\cite{BrzMa:coa}. This is recalled in the following 
\bde\Label{psi2} 
Let $(A,C,\psi)$ be an entwining structure, $\kappa:A\ra k$  
an algebra homomorphism (character), and 
$$
I_\kappa := {\rm span} \{((\kappa\otimes C)\circ\psi)(c\otimes a) -
c\kappa(a) \; | \; a\in A, \; c\in C\}.
$$
Then $B := C/I_\kappa$ is a  coalgebra, and we say that  
$C(B,A,\psi, \kappa)$ is a {\em dual $\psi$-principal  bundle} 
iff the map $cocan_\psi = (C\otimes \kappa\otimes C)\circ (C\otimes
\psi)\circ(\Delta\otimes A): C\otimes A\to C\Box_B C$
is bijective. 
\ede 
Using Theorem~\ref{main2} we can relate dual $\psi$-principal bundles
with $A$-Galois coextensions.
\bpr\Label{cor2} 
For a given entwining structure $(A,C,\psi)$ and an algebra map $\kappa:
A\to k$, the 
 following statements are equivalent:\medskip \\  
\hspace*{2.5mm} 
(1) $C(B,A,\psi,\kappa)$ is a dual $\psi$-principal bundle.\medskip \\  
\hspace*{2.5mm} 
(2) There exists a unique action $\mu_C: C\otimes A\ra C$ such that 
$C\twoheadrightarrow B$ is an $A$-Galois coextension of $B$ by $A$, 
$\psi$ is the canonical entwining 
map,  and $\eps\circ \mu_C = \eps\otimes\kappa$. 
\epr 
\proof (1) $\Rightarrow$ (2). By \cite[Proposition~2.6]{BrzMa:coa}, $C$ is 
a right $A$-module with the action $\mu_C = (\kappa\otimes C)\circ \psi$.
Taking advantage of the second of conditions (\ref{diag.B}) one  verifies that 
$\eps\circ \mu_C = \eps\otimes\kappa$. 
The first of conditions (\ref{diag.B}) implies that 
$C\in \CM_A^C(\psi)$ with structure maps $\Delta$ and $\mu_C$. Using
this fact, explicit form of $\mu_C$,  and (\ref{diag.B})
 one can show that $I_\kappa = I$,
the latter being defined in Lemma~\ref{coa}. Hence $cocan = cocan_\psi$ 
and 
$C\twoheadrightarrow B$ is an $A$-Galois coextension. 
By Theorem~\ref{main2},
$\psi$ is the associated canonical entwining map. Assume now that there
exists another action $\mu'_C$ of $A$ on $C$ satisfying conditions (2).
Then, by Theorem~\ref{main2}, $C\in \CM_A^C(\psi)$ with the structure
maps $\Delta$ and $\mu'_C$. Consequently, since $\eps\circ \mu_C =
\eps\otimes\kappa$, we have
\begin{eqnarray*}
\mu'_C & = & (\eps\otimes C)\circ\Delta\circ \mu'_C \\
& = & (\eps\otimes C)\circ (\mu'_C\otimes C)\circ
(C\otimes\psi)\circ(\Delta\otimes A)\\
& = & (\eps\otimes\kappa\otimes C)\circ(C\otimes\psi)\circ(\Delta\otimes
A)\\
& = & (\kappa\otimes C)\circ\psi = \mu_C.
\end{eqnarray*}

(2) $\Rightarrow$ (1). Since $\psi$ is the canonical entwining map, $C\in 
\CM_A^C(\psi)$ by Theorem~\ref{main2}. The normalisation condition $\eps\circ \mu_C =
\eps\otimes\kappa$ implies that $(\eps\otimes \eps)\circ cocan =
\eps\otimes \kappa$ and, consequently, $\kappa\circ\check{\tau} =
\eps\otimes \eps$. This, in turn, leads to the equality $\mu_C =
(\kappa\otimes C)\circ\psi$. Using this equality, the fact that $C\in 
\CM_A^C(\psi)$, and (\ref{diag.B}) one shows that $I_\kappa = I$.
Hence $cocan_\psi = cocan$ and $cocan_\psi$ is bijective, as required.
\endproof 

\note{the notion of a dual $\psi$-principal bundle 
was introduced in the framework of entwining structures. This 
was possible by noting that the notion of an entwining structure 
has the following self-duality property. 
If one replaces $A$ with $C$, $\Delta$ with $m$, 
$\eta$ with $\eps$ and reverses the order of composition, then formulae 
(\ref{diag.A}) become (\ref{diag.B}) and vice versa. Furthermore, if there is an 
algebra character $\kappa :A\to k$ one can 
define a right action $\mu_C :C\tens A\to C$ by $\mu_C = 
(\kappa\tens C)\circ \psi$. This action is counital 
with respect to $\kappa$. It is also  
left colinear over $B:= C/I_\kappa\,$, where $I_\kappa$ is given 
by (\ref{ikappa}). Indeed,  
\begin{eqnarray*} 
\Delta_\kappa\circ\mu_C &=& 
(\pi_\kappa\tens C)\circ\Delta\circ(\kappa\tens C)\circ\psi 
\\ &=& 
(\pi_\kappa\tens C)\circ(\kappa\tens C\tens C)\circ (A\tens \Delta)\circ\psi 
\\ &=& 
(\pi_\kappa\tens C)\circ(\kappa\tens C\tens C)\circ(\psi\tens 
C)\circ(C\tens\psi)\circ(\Delta\tens A) 
\\ &=& 
(\pi_\kappa\circ\mu_C\tens C)\circ(C\tens\psi)\circ(\Delta\tens A) 
\\ &=& 
 (\pi_\kappa\tens\kappa\tens A)\circ(C\tens\psi)\circ(\Delta\tens A) 
\\ &=& 
(\pi_\kappa\tens\mu_C)\circ(\Delta\tens A)  
\\ &=& 
(B\otimes\mu_C)\circ(\Delta_\kappa\tens A). 
\end{eqnarray*} 
We used the definition of $\mu_C$ to derive the first, fourth and sixth 
equalities, the definition of $I_\kappa$ to derive the fifth one, and the 
first of equations (\ref{diag.B}) to establish the third equality. 
By Lemma~\ref{colem} $B$ is a coalgebra, and we can consider 
algebra-Galois coextensions whose right action map $\mu_C$ is determined 
by $\psi$ as described above. }
 
\note{\bde\Label{copsi} 
Let $(A,C,\psi)$ be an entwining structure, $\kappa :A\to k$  an  
algebra character, and $\mu_C:=(\kappa\ot C)\ci\psi$ a right action 
of $A$ on $C$. We say that  
$C\twoheadrightarrow B$ is a {\em $\psi$-Galois coextension} of 
$B=C/I_\kappa$ by an algebra $A$ iff $C\twoheadrightarrow B$ is an $A$-Galois 
coextension, i.e., 
iff the canonical morphism $cocan  : 
C\otimes A\to C\Box_BC$ in $^C\!\CM_A$ is an isomorphism. 
(We call the underlying quantum-space structure of $(A,B,C,\psi,\kappa)$  
 a dual $\psi$-principal bundle.) 
\note{We say that 
 $(C,B,A,\psi,\kappa)$ is a {\em dual $\psi$-principal 
 bundle} over $B$ or, equivalently, that $C\twoheadrightarrow B$ is an 
 {\em $A$-Galois coextension} if $C\tens A \cong C\Box_{B} C$ 
 as objects in $^C\!\CM_A$  by the canonical map 
 $\can' :C\otimes A\to C\Box_{B} C$, $\can':  c\tens a \mapsto c\sw 
 1\tens \mu_C(c\sw 2, a)$.} 
\ede }
 
%As in the previous section, we have 
 
\note{\parindent=5pt 
\bco 
For an algebra $A$, a coalgebra $C$ and an algebra 
character $\kappa :A\to k$, the following statements are equivalent: }
 
%\hspace*{2.5mm} (1) $(A,B,C,\psi,\kappa)$ is a $\psi$-Galois coextension. 
 
\note{\hspace*{2.5mm} (2) $C\twoheadrightarrow B$ is an $A$-Galois coextension of  
$B:=C/I_\kappa$ by $A$. 
\eco }

\section{Appendix} 
 
We say that  two subgroups $G_1$ and $G_2$ of a group $G$ generate $G$ 
if any element of $G$ can be written as a finite length word whose 
letters are elements of $G_1$ or $G_2\,$. The (dual) coalgebra version 
of this concept is given in the following definition: 
 
\bde\Label{codef} 
Let $C$ be a coalgebra and $I_1\,$, $I_2$ its coideals. Let 
$\wp_{(i)}$ denote the composite map 
\[ 
C\st{\Delta_n}{-\!\!\!\lra}C^{\ot n+1}\st{\pi_{i_1}\ot\dots\ot\pi_{i_{n+1}}} 
{-\!\!\!-\!\!\!-\!\!\!-\!\!\!-\!\!\!\lra}C/I_{i_1}\ot\dots\ot C/I_{i_{n+1}}
\, , 
\] 
where $\hD_n(c):=c\1\ot\dots\ot c_{(n+1)}$,  
$(i):=(i_1,\dots,i_n)\in\{1,2\}^{\times n}$ is a finite multi-index,  
and each $\pi_{i_k}$ is a canonical surjection. 
We say that the quotient coalgebras $C/I_1$ and $C/I_2$ {\em cogenerate} $C$ 
iff $\bigcap_{(i)\in{\cal M}_f}\mbox{\em Ker}\,\wp_{(i)}=0$, 
where ${\cal M}_f$ 
is the space of all finite multi-indices.  
We write then  $(C/I_1)\cdot(C/I_2)=C$. 
\ede 
\parindent=5mm 
 
Observe that the above construction is closely related to the wedge 
construction of~\cite{swe}. 
In the group situation it is clear that what is invariant under both 
generating subgroups is invariant under the whole group, and vice-versa. 
Below is the (dual) coalgebra version of this classical phenomenon. 
 
\bpr\Label{coin} 
Let $C$ be a coalgebra, $I_1$ and $I_2$ 
coideals of $C$, and $A$ a right $C$-comodule. Then, defining coinvariants
as in Definition~\ref{coinv}, we have: 
\[ 
(C/I_1)\cdot(C/I_2)=C\;\;\imp\;\; A\sp{co C}
=A\sp{co (C/I_1)}\cap A\sp{co 
(C/I_2)}. 
\] 
\epr 
\proof Clearly, we always have 
$A\sp{co C}\inc A\sp{co (C/I_1)}\cap A\sp{co (C/I_2)}$. Assume now that 
there exists $b\in A\sp{co (C/I_1)}\cap A\sp{co (C/I_2)}$ such that 
$b\not\in A\sp{co C}$. Then there also exists $a\in A$ such that 
$0\not=\hD_A(ba)-b\hD_A(a)=:\sum_{j\in{\cal J}}f_j\ot h_j$, where 
$\{f\sb\alpha\}\sb{\alpha\in{\cal A}}$ is a basis of $A$ and 
$\{h_j\}\sb{j\in{\cal J}}$ is a non-empty set which does not contain zero. 
Furthermore, for any $(i)\in{\cal M}_f$, we have: 
\bea 
&& 
(A\ot\wp_{(i)})(\mbox{$\sum_{j\in{\cal J}}$}f_j\ot h_j) 
\nonumber \\ && = 
\llp(A\ot\pi_{i_1}\ot\cdots\ot\pi_{i_n})\circ(A\ot\hD\sb{n-1})\lrp 
((ba)\0\ot(ba)\1-ba\0\ot a\1) 
\nonumber\\ && = 
%\llp(A\ot\pi_{i_1}\ot\cdots\ot\pi_{i_n})
%\circ(A\ot\hD\sb{n-1})\circ\hD_A\lrp(a)-a\ot e\sp{\otimes n} 
%\nonumber\\ && = 
(ba)\0\ot\pi_{i_1}((ba)\1)\ot\cdots\ot\pi_{i_n}((ba)_{(n)})
-ba\0\ot\pi_{i_1}(a\1)\ot\cdots\ot\pi_{i_n}(a_{(n)})
\nonumber\\ && = 
%\llp(A\te\pi_{i_1})\cc\hD_A\ot C\sp{\otimes (n-1)}\lrp 
%\llp a\0\ot\pi_{i_2}(a\1)\ot\cdots\ot\pi_{i_n}(a_{(n-1)})\lrp -a\ot 
%e\sp{\otimes n} 
%\nonumber\\ && = 
\llp[(A\te\pi_{i_1})\cc\hD_A\ot C\sp{\otimes (n-1)}]\ci\cdots\ci 
(A\ot\pi_{i_n}) 
\ci\hD_A\lrp(ba)
\nonumber\\ && \phantom{=} -
ba\0\ot\pi_{i_1}(a\1)\ot\cdots\ot\pi_{i_n}(a_{(n)})
\nonumber\\ && = 
0\nonumber 
\eea 
Consequently, by the linear independence of $f_j,\, j\in{\cal J}$, we have 
$\wp\sb{(i)}(h_j)=0$ for $j\in{\cal J}$. Hence, as this is true for any 
$(i)\in{\cal M}_f$, we obtain 
$\bigcap_{(i)\in{\cal M}_f}\mbox{Ker}\wp\sb{(i)}\not=0$, 
as needed. 
\epf 
 
\begin{center} 
{ACKNOWLEDGEMENTS}\vspace{-16pt} 
\end{center} 
This work was supported by the EPSRC grant GR/K02244, the University of 
 \L\'od\'z grant 505/582 and the Lloyd's Tercentenary Foundation  (T.B.), and 
by  the  NATO postdoctoral fellowship and KBN grant \mbox{2 P301 020 07} 
(P.M.H.). The authors are very grateful to Mitsuhiro Takeuchi for providing
a definition of coinvariants that does not require the
existence of a group-like element, and his suggestion to use it in this 
paper.

\end{document}